# SUBFACTORS OF INDEX LESS THAN 5, PART 3: QUADRUPLE POINTS

MASAKI IZUMI, VAUGHAN F. R. JONES, SCOTT MORRISON, AND NOAH SNYDER

Abstract. One major obstacle in extending the classification of small index subfactors beyond $3+\sqrt{3}$ is the appearance of infinite families of candidate principal graphs with 4-valent vertices (in particular, the "weeds" $\mathcal{Q}$ and $\mathcal{Q}'$ from Part 1 [MS10]). Thus instead of using triple point obstructions to eliminate candidate graphs, we need to develop new quadruple point obstructions. In this paper we prove two quadruple point obstructions. The first uses quadratic tangles techniques and eliminates the weed $\mathcal{Q}'$ immediately. The second uses connections, and when combined with an additional number theoretic argument it eliminates both weeds $\mathcal{Q}$ and $\mathcal{Q}'$. Finally, we prove the uniqueness (up to taking duals) of the 3311 Goodman-de la Harpe-Jones subfactor using a combination of planar algebra techniques and connections.

## Contents

## 1. Introduction

This is the third paper in a series of four papers (along with [MS10, MPPS10, PT10]) in which we extend the previously known classification of subfactors of index less than $3+\sqrt{3}$ [Haa94, Bis98, AH99, AY09, BMPS09] up to index 5.

**Theorem 1.1.** *There are exactly ten subfactor planar algebras other than Temperley-Lieb with index between 4 and 5: the Haagerup planar algebra and its dual [AH99], the extended Haagerup planar algebra and its dual [BMPS09], the Asaeda-Haagerup planar algebra [AH99] and its dual, the 3311 Goodman-de la Harpe-Jones planar algebra [GdlHJ89] and its dual, and Izumi's self-dual 2221 planar algebra [Izu01] and its complex conjugate.*

In the first paper of this series [MS10], we gave an initial classification result which says that any exceptions to Theorem 1.1 must lie in one of finitely many families. To complete the classification we rule out the remaining families in the subsequent papers. In this paper we eliminate two families with 4-valent vertices.

In order to state our main results we rapidly recall some terminology from Part 1. A translation of a graph pair is used to indicate a graph pair obtained by increasing the supertransitivity by an even integer (the supertransitivity is the number of edges between the initial vertex and the first vertex of degree more than two). An extension of a graph pair is a graph pair obtained by extending the graphs in any way at greater depths (i.e. adding vertices and edges at the right), perhaps even infinitely.

Combining the results of the first two papers we have

**Theorem 1.2** (From [MS10, MPPS10]). *The principal graph of any subfactor of index between* 4 *and* 5 *is a translate of one of an explicit finite list of graph pairs, which we call the* vines*, or is a translated extension of one of the following graph pairs, which we call the* weeds*.*

$\mathcal{Q} =$ 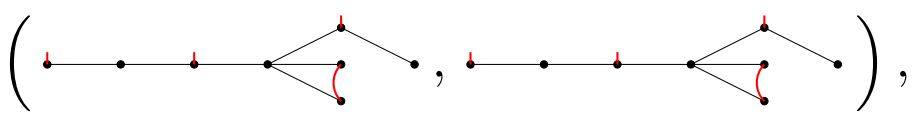

$\mathcal{Q}' =$ 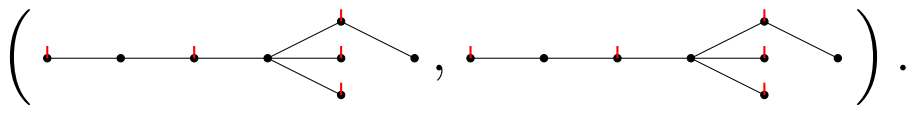

(As in [MS10], the trivial bimodule always appears as the leftmost vertex of a principal graph. Dual pairs of bimodules are indicated in red at even depths, and by matching up vertices on the two graphs at corresponding heights at odd depths.)

Since $\mathcal{Q}$ and $\mathcal{Q}'$ differ only in their dual data, we refer to the underlying graph, forgetting dual data, by

$\mathcal{Q}^\bullet =$ 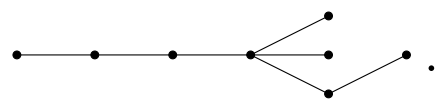

The goal of this paper is to eliminate the two remaining weeds $\mathcal{Q}$ and $\mathcal{Q}'$. We do this by proving three theorems.

**Theorem 1.3.** *There is no subfactor (of any index) whose principal graph is a translated extension of* $\mathcal{Q}$*.*

We prove Theorem 1.3 using planar algebra techniques based on [Jon01, Jon03, BMPS09]. In particular, we show that any subfactor planar algebra whose principal graph is a translated extension of $\mathcal{Q}$ has a proper planar subalgebra which begins with a triple point at the same depth. This quickly leads to a contradiction using Ocneanu's triple point obstruction [Haa94].

**Theorem 1.4.** *Any subfactor (of any index) whose principal graph is a translated extension of* $\mathcal{Q}^\bullet$ *has principal graph* 3311 *and index* $3 + \sqrt{3}$*.*

We prove this result in two steps. First we show, using connections, that any such principal graph must have the same index as the $nn11$ graph and be finite depth. Then, using number theoretic techniques from [CMS11], we prove that the graph norm of the $nn11$ graph is non-cyclotomic and thus there can be no finite depth subfactor with that index. Note that this theorem is strictly stronger than the previous theorem, but nonetheless we have included both arguments since they illustrate very different approaches to quadruple point obstructions.

**Theorem 1.5.** *There are exactly two subfactors with principal graph* 3311*, namely the GHJ subfactor and its dual.*

This result is originally due to Rehren in unpublished notes from 1994, who proved directly that of the 1-parameter family of connections only two of them can be flat. We give an independent proof which combines connections arguments with planar algebraic techniques. In particular, instead of checking flatness directly we show how to read off the rotational eigenvalues from the connection and exploit that to determine which connections could be flat. This relationship between rotational eigenvalues and connections should be applicable more generally. In order to state this result, we first need a characterization of a particular gauge choice for the connection, resulting the the following definition of a 'diagrammatic branch matrix'.

**Definition 1.6.** Given an $n$-supertransitive principal graph pair $\Gamma$, with a $k{+}2$-valent initial branch point, with no multiple edges, and a flat biunitary connection on $\Gamma$, the branch matrix is the $(k+2)$-by-$(k+2)$ unitary matrix consisting of the values of the connection on loops between the initial branch points on the two principal graphs.

A diagrammatic branch matrix is a branch matrix for which the top left entry of $U$ is a real number with sign $(-1)^{n+1}$, and the other entries in the first row and first column are positive real numbers.

In fact any flat biunitary connection has a unique gauge group representative such that its branch matrix is a diagrammatic branch matrix.

**Theorem 1.7.** *Suppose $U$ is the diagrammatic branch matrix for an $n$-supertransitive principal graph pair $\Gamma$, having a $k+2$-valent initial branch point with no multiple edges, and a flat biunitary connection. Let $\rho$ be the 'two-click' rotation $P_{n+1,+} \to P_{n+1,+}$ in the corresponding planar algebra $P$.*

*The eigenvalues for $U^t U$ with multiplicity are: two 1's together with the eigenvalues for $\rho$ acting on the perpendicular complement of Temperley-Lieb in the $n$-box space. In particular, $tr(U^t U) - 2$ is the sum of the eigenvalues of $\rho$.*

This result is proved in §5.2.

We recall that Part 4 of this series, [PT10], deals with the remaining vines mentioned in Theorem 1.2 above, and thus along with the results contained here completes the proof of Theorem 1.1.

We will assume that readers are familiar with planar algebras [Jon99] and with connections [Ocn88, EK98]; see the background section to [MPPS10] and [BMPS09, §2.1] for more details. In Section 2 we also assume familiarity with the main ideas of [Jon03] and in Section 4 we assume familarity with the main ideas of [CMS11].

We would like to thank Stephen Bigelow, David Penneys, Emily Peters, and Karl-Henning Rehren for helpful conversations. Noah and Scott would like to thank Kyoto University where some of this research was done. During this work, Scott Morrison was at the Miller Institute for Basic Research at UC Berkeley, and Noah Snyder held an NSF Postdoctoral Fellowship at Columbia University. Masaki Izumi is supported by the Grant-in-Aid for Scientific Research (B) 22340032, JSPS. Vaughan Jones is supported by DMS 0856316 and the NZIMA. We would also like to acknowledge support from the DARPA HR0011-11-1-0001 grant.

## 2. Quadratic tangles

The connections approach to the standard invariant of a subfactor emphasizes the fundamental role of projections. A key observation from [Jon01, Jon03] is that in the planar algebraic description [Jon99] of the standard invariant you can instead build the theory around the eigenvectors of the rotation operator $\rho$. Since rotational eigenvalues are necessarily roots of unity this gives obstructions that aren't readily visible in the connections approach. Throughout this paper $\rho$ refers to the 'two-click' rotation $\rho : P_{n+1} \to P_{n+1}$, when $P$ is an $n$-supertransitive subfactor planar algebra. The key idea in this section is to show that a subfactor planar algebra whose principal graph is a translated extension of $\mathcal{Q}$ has a rotational eigenvector which behaves like it comes from a triple-point. Thus any such subfactor principal graph can be eliminated using a triple point obstruction.

We give two slightly different proofs of Theorem 1.3 following this general outline. The first shows that the proof of the triple-point obstruction from [Jon03] can be applied to one of the rotational eigenvectors. The second proof is slightly more involved, but conceptually explains why you can use a triple point obstruction. Using a variant of the jellyfish algorithm [BMPS09], we show that one of the rotational eigenvectors generates a proper planar subalgebra and the principal graph of this planar subalgebra begins with a triple point. In this version, we can get a contradiction using Ocneanu's original triple point obstruction [Haa94]. The key point that makes both of these arguments work is that we can exploit the fact that the projections one past the split in $\mathcal{Q}$ are not self-dual to gain control over the change of basis between the projection basis and the rotational eigenvalue basis.

Let $n$ be the depth of the quadruple point. Consider the projections $p$, $q$ and $r$ in the unshaded $n+1$ box space, with $p$ and $q$ corresponding to the univalent vertices of the principal graph (and so $p^* = q$). The perpendicular complement of Temperley-Lieb in the $n+1$-box space is 2-dimensional. Since Temperley-Lieb is preserved by rotation (and rotation is unitary) we can diagonalize the rotation on the perpendicular complement giving two perpendicular rotational eigenvectors $S$ and $T$. By [Jon03, Theorem 5.2.3] the rotational eigenvalues are distinct, and thus $S$ and $T$ are well-defined up to rescaling. Let $f$ denote the $n+1$-strand Jones-Wenzl idempotent. Similarly in the shaded $n+1$ box spaces we have $p', q', r', f', S', T'$. We will not assume that $S'$ is exactly $\rho^{\frac{1}{2}}(S)$ (but it must be a multiple because both are eigenvectors with the same eigenvalue), as it is more convenient to be able to normalize $S$ and $S'$ independently.

**Lemma 2.1.** *Without loss of generality, $\pi$-rotation acts on $S$ by $+1$ and on $T$ by $-1$. Furthermore, the rotational eigenvectors can be written in terms of the projections as $T = p-q$ and $S = -(1-2a)p-(1-2a)q+2ar$, where $a = tr(p) = tr(q) = [n+1]/[2]$. The same results hold in the other shading. Finally, $\rho^{\frac{1}{2}}(S) = \omega^{-\frac{1}{2}}S'$ where $\omega$ is the rotational eigenvalue of $S$.*

*Proof.* Inside Temperley-Lieb is the span of diagrams with smaller waists (that is where at least one strand from the bottom connects again to the bottom). Note that $f$, $p$, $q$, $r$, $S$, and $T$ all lie in the perpendicular complement to this space. In particular we have the following change of basis, for some values of $b, b', c$ and $c'$ (the value of $a$ is determined by the traces).

$$\begin{aligned} p &= af + bS + cT \\ q &= af + b'S + c'T \\ r &= (1-2a)f - (b+b')S - (c+c')T \end{aligned}$$

Now $\pi$-rotation acts on $S$ and $T$ by signs. Since $\pi$-rotation interchanges $p$ and $q$ we have that $b' = \pm b$ and $c' = \pm c$. At least one of these signs must be negative, because $p$ and $q$ are distinct. Furthermore, they can not both be negative because $r$ is not a multiple of $f$. Without loss of generality, let $S$ have $\pi$-rotational eigenvalue $+1$ and $T$ have $\pi$-rotational eigenvalue $-1$. Furthermore, normalize $S$ and $T$ so that $c = b = 1/2$. Thus,

$$\begin{aligned} p &= af + \frac{1}{2}S + \frac{1}{2}T \\ q &= af + \frac{1}{2}S - \frac{1}{2}T \\ r &= (1-2a)f - S \end{aligned}$$

Inverting, we get the change of basis the other direction.

$$\begin{aligned} f &= p + q + r \\ T &= p - q \\ S &= -(1-2a)p - (1-2a)q + 2ar \end{aligned}$$

For the last statement note that since $S$ and $S'$ have the same eigenvalue for $\pi$-rotation, we have that $\rho^{\frac{1}{2}}(S)$ is some multiple of $S'$. To work out the actual multiple note that by sphericality $\mathrm{tr}(\rho^{\frac{1}{2}}(S)^2) = \omega^{-1}\mathrm{tr}(S^2)$. (Note that there's a sign ambiguity depending on the choice of square root of $\omega$. Either choice works and the graph automorphism interchanging $p'$ and $q'$ changes the choice.) □

It follows that $S^2$ lies in the span of $f$ and $S$. This is the first hint that $S$ behaves like part of a triple-point planar algebra.

Recall from [Jon03] that if $A$ and $B$ are two $k$-boxes we can build a $k+1$-box $A \circ B$ by connecting only $k-1$ of the strands from $A$ to $B$.

**Lemma 2.2.** *$S \circ S$ is perpendicular to the space of annular consequences of $T$, and $S' \circ S'$ is perpendicular to the space of annular consequences of $T'$.*

*Proof.* It is immediate that $S \circ S$ is perpendicular to all but two of the annular consequences of $T$, as one of the $S$'s gets capped off. For the remaining two annular consequences of $T$, the inner products are $\mathrm{tr}(S^2 T)$ and $\mathrm{tr}(\rho^{\frac{1}{2}}(S)^2 \rho^{\frac{1}{2}}(T))$. Up to a non-zero scalar the latter is the same as $\mathrm{tr}(S'^2 T')$. In either case we compute: $\mathrm{tr}(S^2 T) = \mathrm{tr}(((1-2a)(p+q)+4a^2 f)(p-q)) = \mathrm{tr}((1-2a+4a^2)(p-q)) = 0$. □

The above two lemmas shows that $S$ "decouples" from $T$ in some sense. At this point there are two ways to complete the argument. We rapidly sketch the first argument and provide a bit more detail for the second.

*First proof of Theorem 1.3.* We apply the argument from [Jon03, Theorem 5.1.11] to $S$ mutatis mutandis. (Specifically at a few points we use the fact that $T$ and $T'$ do not appear in $S^2$, $S'^2$, $S \circ S$, or $S' \circ S'$.) The numbers playing the role of $r$ and $\check{r}$ in [Jon03, Theorem 5.1.11] are $\frac{\dim p + \dim q}{\dim r}$ (or its inverse) and $\frac{\dim p' + \dim q'}{\dim r'}$ (or its inverse). It's easy to see that none of these ratios can be $\frac{[n+2]}{[n]}$, which gives a contradiction. □

Now we turn to the second proof. Since the first argument uses a triple point obstruction, it would be more satisfying to have an argument which works via a subfactor planar algebra with an initial triple point. In the second proof we show that $S$ generates such a planar subalgebra.

**Lemma 2.3.** *The following annular consequence of $S$ is a linear combination of Temperley-Lieb diagrams, other annular consequences of $S$, and $S \circ S$.*

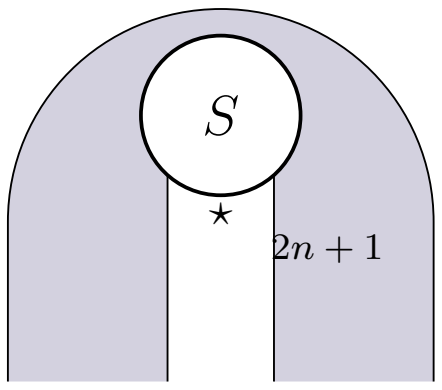

*The analogous result also holds for $S'$.*

*Proof.* By counting dimensions (and using the fact that annular consequences are linearly independent from Temperley-Lieb and other annular consequences [Jon01]) we see that Temperley-Lieb diagrams together with annular consequences of $S$ and $T$ is a basis for the $n+2$ box space. We write $S \circ S$ in this basis. Since $S \circ S$ is perpendicular to annular consequences of $T$ this lies in the span of Temperley-Lieb and annular consequences of $S$. Thus, we need only show that the coefficient of this specific annular consequence in $S \circ S$ is nonzero. This is just an inner product calculation. Note that the following element is perpendicular to every other basis element.

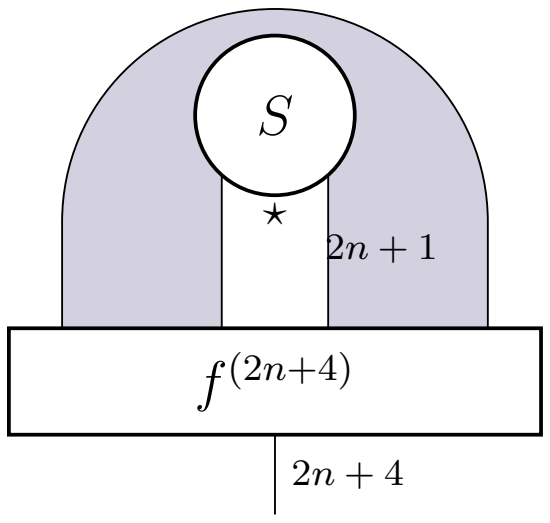

Thus we need only evaluate the following diagram and see that it is nonzero.

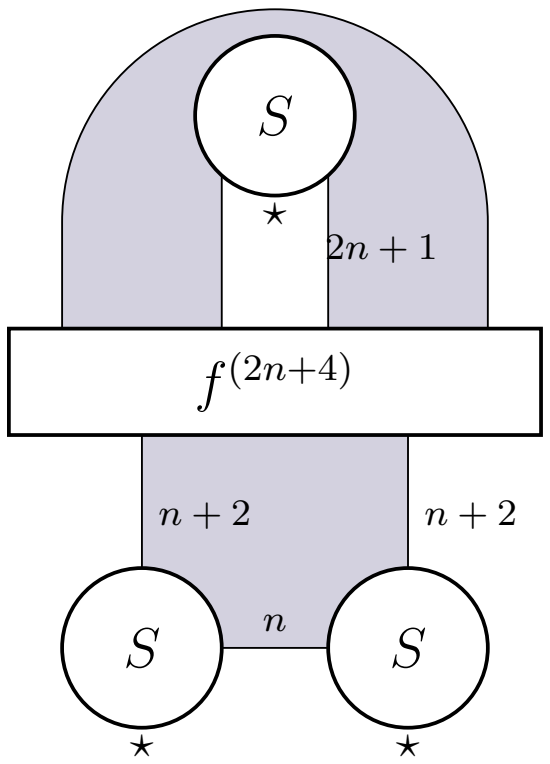

By the calculation in [BMPS09, Lemma 5.5] we see that this is

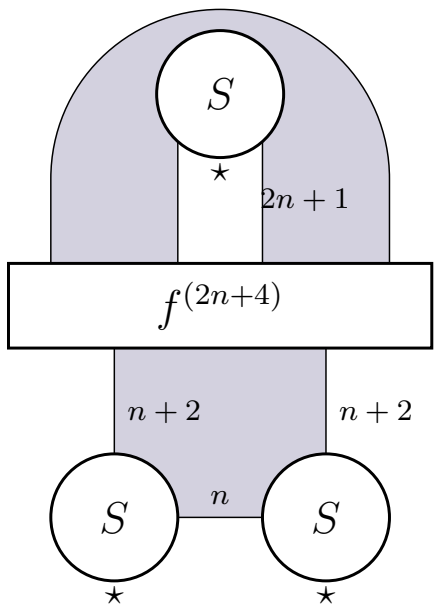

$$= \mathrm{tr}(\rho^{\frac{1}{2}}(S)^3) + \mathrm{tr}(S^3) + \omega^{-1}\mathrm{tr}(S^3)$$

$$= \omega^{-3/2}\mathrm{tr}(S'^3) + \mathrm{tr}(S^3) + \omega^{-1}\mathrm{tr}(S^3)$$
$$= (\omega^{-\frac{3}{2}} + 1 + \omega^{-1})\mathrm{tr}(S^3) \neq 0.$$

□

**Lemma 2.4.** *The element $S$ generates a planar subalgebra $P(S)$ whose $n+1$ box space is spanned by Temperley-Lieb and $S$, and whose $n+2$ box space is spanned by Temperley-Lieb and the annular consequences of $S$.*

*Proof.* This follows from the jellyfish algorithm [BMPS09, §4] combined with the argument in [Big10, Theorem 5.1]. In fact, using Bigelow's ideas you can find an explicit spanning set (of diagrams in "jellyfish form") for all every box space in $P(S)$.

Start with an arbitrary planar composition of $S$'s. First use the relations proved in the last lemma to pull all the $S$'s to the top of the diagram (along the way the number of $S$'s will increase). If two $S$'s are connected by $n+1$ strands, then simplify the diagram using the formula for $S^2$ or for $S'^2$ (after possibly repeatedly applying the rotational eigenvector relation). Hence any diagram in $P(S)$ is in the space of diagrams of jellyfish form: with all $S$'s at the top of the diagram and no two $S$'s connected by $n+1$ strands. (Thus far we've mimicked [BMPS09, §4].)

Since the $S$'s are all at the top of the diagram, we can think of them as lying on a line with each $S$ having well-defined neighbor on the left and on the right. Following [Big10, Lemma 4.3 and Theorem 5.1], we claim that any diagram in jellyfish form has each $S$ only attached to its neighbors and the outside. If there were two non-adjacent $S$'s connected by a strand, then somewhere in between them there would be an $S$ which was only attached to its neighbors and not to the outside. Such an $S$ would have to be connected to one of its neighbors by $n+1$ strands.

Any diagram in jellyfish form with two $S$'s must have at least $2n+4$ boundary strands, because the leftmost $S$ must have $n+2$ boundary strands and the rightmost $S$ must also have $n+2$ boundary strands. Similarly, any diagram in jellyfish form with three or more $S$'s must have at least $2n+6$ boundary strands (the leftmost and rightmost each contribute $n+2$, while each $S$ in the middle must contribute 1, and any diagram has an even number of boundary points). Thus the only diagrams in jellyfish form with $2n+2$ boundary points are in Temperley-Lieb or are a single $S$, while the only diagrams in jellyfish form with $2n+4$ boundary points are in Temperley-Lieb, annular consequences of $S$, or $S \circ S$. □

We can now apply triple point obstructions to $P(S)$.

*Second proof of Theorem 1.3.* Consider the planar algebra $P(S)$ generated by $S$. We claim that both the principal and dual principal graphs of $P(S)$ are extensions of translations of 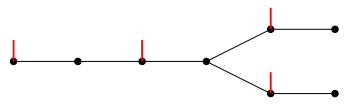. By the previous lemma they must be a extensions of translations of 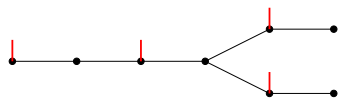 or of 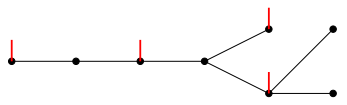. However, the simple projections at depth $n+1$ are $p+q$ and $r$. Now, in $P(S)$, the two idempotents are $p+q$ and $r$. The products $(p+q)X$ and $rX$ each have at least two summands, so there are no univalent vertices at depth $n+1$. Finally either the triple point obstruction from [Haa94] or from [Jon03] completes the proof. □

## 3. Connections

In this section, we study biunitary connections on certain principal graphs with quadruple points. In particular, we show in Corollary 3.3 that any extension of the $(n-3)$ translate of $\mathcal{Q}^\bullet$ has the same index as the $nn11$ graph. This implies that there are no infinite depth extensions, because any such an extension would have index strictly greater than that of $nn11$. In addition we describe all connections on 3311.

As usual in the connections business, we need some linear algebra lemmas about unitary matrices. For triple-point obstructions the key elementary lemma is the triangle inequality. In particular, if you have two orthogonal 3-dimensional vectors whose components have known sizes but unknown phases, then the triangle inequality can give a contradiction. For quadruple-point obstructions the key observation is that rhombuses are parallelograms. Thus if you have two orthogonal 4-dimensional vectors whose components all have the same size but unknown phases, then you get a relation between the phases.

**Lemma 3.1.** *Let $A, B, C, D$ and $\beta$ be positive numbers satisfying*

$$\beta A = B + 2C + D,\ \beta C = A,\ A^2 - \beta AB + B^2 > 0$$
$$\beta > 2,\ \beta B > A,\ \beta D > A.$$

*Let $\xi_i, \eta_i$, for $i = 1, \ldots, 4$, be unit complex numbers, and $x$ and $y$ be complex numbers. If the matrix*

$$\frac{1}{A}\begin{pmatrix} y & \sqrt{BC} & \sqrt{BC} & \sqrt{BD} \\ \sqrt{BC} & C\xi_1 & C\xi_2 & \sqrt{CD}\eta_1 \\ \sqrt{BC} & C\xi_3 & C\xi_4 & \sqrt{CD}\eta_2 \\ \sqrt{BD} & \sqrt{CD}\eta_3 & \sqrt{CD}\eta_4 & x \end{pmatrix}$$

*is unitary, then $x = -\bar{y}$ and $B = D$.*

*Proof.* From the orthogonality of the first and last rows and that of the first and last columns, we obtain

$$\bar{y} + C(\eta_1 + \eta_2) + x = \bar{y} + C(\eta_3 + \eta_4) + x = 0,$$

which implies that $\eta_1 + \eta_2 = \eta_3 + \eta_4$, and so $\Re(\eta_2\bar{\eta_1}) = \Re(\eta_4\bar{\eta_3})$. Without loss of generality (since we can switch the second and third rows) we may assume $\eta_2\bar{\eta_1} = \eta_4\bar{\eta_3}$, and we denote this quantity $\alpha$. Then we have $(1+\alpha)(\eta_1 - \eta_3) = 0$. If $\alpha = -1$, we have $x = -\bar{y}$. Considering the norm of first and last row, we see

$$|x|^2 = A^2 - \beta AD + D^2 \tag{3.1}$$

$$|y|^2 = A^2 - \beta AB + B^2 \tag{3.2}$$

and thus $B(B - \beta A) = D(D - \beta A)$. This is only satisfied if $B + D = \beta A$, which is impossible as it implies $C = 0$, or $B = D$, as desired.

Now assume $\alpha \neq -1$. Then we have $\eta_3 = \eta_1$, which we will just denote by $\eta$ from now on. By the orthogonality of the first and second rows and that of the first and second columns, we obtain

$$\bar{y} + C(\xi_1 + \xi_2) + D\eta = \bar{y} + C(\xi_1 + \xi_3) + D\eta = 0$$

and so $\xi_2 = \xi_3$. Now our unitary matrix is

$$\frac{1}{A}\begin{pmatrix} y & \sqrt{BC} & \sqrt{BC} & \sqrt{BD} \\ \sqrt{BC} & C\xi_1 & C\xi_2 & \sqrt{CD}\eta \\ \sqrt{BC} & C\xi_2 & C\xi_4 & \sqrt{CD}\alpha\eta \\ \sqrt{BD} & \sqrt{CD}\eta & \sqrt{CD}\alpha\eta & x \end{pmatrix}$$

whose unitarity is equivalent to

$$\bar{y} + C(\xi_1 + \xi_2) + D\eta = 0, \tag{3.3}$$

$$\bar{y} + C(\xi_2 + \xi_4) + D\alpha\eta = 0, \tag{3.4}$$

$$\bar{y} + C(1+\alpha)\eta + x = 0, \tag{3.5}$$

$$B + C(\xi_2\bar{\xi_1} + \bar{\xi_2}\xi_4) + D\alpha = 0, \tag{3.6}$$

$$B + C\bar{\eta}(\xi_1 + \bar{\alpha}\xi_2) + \eta\bar{x} = 0, \tag{3.7}$$

$$B + C\bar{\eta}(\xi_2 + \bar{\alpha}\xi_4) + \alpha\eta\bar{x} = 0 \tag{3.8}$$

and Equations (3.1) and (3.2) above. Suppose that $\alpha = 1$, then Equations (3.3) and (3.4) imply that $\xi_1 = \xi_4$, so Equation (3.6) implies that

$$\xi_1\bar{\xi_2} + \bar{\xi_1}\xi_2 = -\frac{B+D}{C} = -\frac{\beta A - 2C}{C} = 2 - \beta^2.$$

This contradicts $\beta > 2$. Now we know that $\alpha \neq 1$, and Equations (3.3), (3.4), (3.7) and (3.8) imply

$$C(\xi_1 - \xi_4) + D(1-\alpha)\eta = 0,$$
$$C(\xi_1 - \bar{\alpha}^2\xi_4) + B(1-\bar{\alpha})\eta = 0.$$

Solve these, we obtain

$$\xi_1 = -\frac{B + \bar{\alpha}D}{C(1+\bar{\alpha})}\eta,$$
$$\xi_4 = -\frac{B + \alpha D}{C(1+\bar{\alpha})}\eta.$$

Since $|\xi_1| = |\xi_4| = 1$, we obtain

$$B^2 + D^2 + BD(\alpha + \bar{\alpha}) = C^2(2 + \alpha + \bar{\alpha}).$$

As we have $\alpha \neq \pm 1$ and $BD - C^2 = (\beta B\beta D - A^2)/\beta^2 > 0$, we obtain

$$B^2 + D^2 - 2C^2 = (C^2 - BD)(\alpha + \bar{\alpha}) < 2(BD - C^2),$$

and $(B-D)^2 < 0$, which is a contradiction. □

**Corollary 3.2.** *If the principal and dual principal graphs of a subfactor are each translated extensions of $Q^\bullet$, then the two bivalent vertices adjacent to the quadruple point must have the same dimension.*

*Proof.* The existence of a subfactor implies the existence of a biunitary connection. This biunitary connection assigns to the pair of initial quadruple points a certain 4-by-4 matrix. The normalization property of biunitary connections tells you that the entries of this matrix have the same magnitudes as the entries of the matrix in the above lemma where $\beta$ is the square root of the index, $A$ is the dimension of the bimodule at the quadruple point, $B$ and $D$ are the dimensions of the objects at the 2-valent vertices, and $C$ is the dimension of the objects at the 1-valent vertices. After applying a gauge transformation we can assume that the first row and column of the matrix are real. □

**Corollary 3.3.** *For any $(n-3)$-translated extension of $\mathcal{Q}^\bullet$, the index is the same as the index of the $nn11$ graph.*

*Proof.* A positive dimension function for the principal graph induces an FP eigenvector for the $nn11$ graph by assigning the same dimension to all vertices up to depth $n+1$, and then the rest by symmetry. □

We now turn our attention to 3311.

**Lemma 3.4.** *Let $A = (1+\sqrt{3})\sqrt{3+\sqrt{3}}$, $B = 2+\sqrt{3}$, $C = 1+\sqrt{3}$ and let $\xi_i$, $\eta_i$ be unit complex numbers. The matrix*

$$\frac{1}{A}\begin{pmatrix} 1 & \sqrt{BC} & \sqrt{BC} & B \\ \sqrt{BC} & C\xi_1 & C\xi_2 & \sqrt{BC}\eta_1 \\ \sqrt{BC} & C\xi_3 & C\xi_4 & -\sqrt{BC}\eta_1 \\ B & \sqrt{BC}\eta_3 & -\sqrt{BC}\eta_3 & -1 \end{pmatrix}$$

*is unitary if and only if*

$$\begin{aligned} \xi_1 &= \frac{-(1+\bar{\alpha}) + (-1+\bar{\alpha})\eta_1 B}{2C} \\ \xi_2 &= \frac{(-1+\bar{\alpha}) - (1+\bar{\alpha})\eta_1 B}{2C} \\ \xi_3 &= \frac{-(1+\alpha) + (1-\alpha)\eta_1 B}{2C} \\ \xi_4 &= \frac{(-1+\alpha) + (1+\alpha)\eta_1 B}{2C} \\ \eta_3 &= \frac{1+\eta_1 B}{\eta_1 + B}\alpha \end{aligned}$$

*where*

$$\alpha = \pm i \frac{2\sqrt{3} + \eta_1 - \bar{\eta_1}}{\left|2\sqrt{3} + \eta_1 - \bar{\eta_1}\right|}.$$

*Remark.* The matrix in the hypotheses of this lemma is exactly the matrix of Lemma 3.1 after substituting $y = 1$ and observing that unitarity implies $\eta_2 = -\eta_1$ and $\eta_4 = -\eta_3$

*Proof.* From the orthogonality of the second and third rows, we have

$$\bar{\xi_1}\xi_3 + \bar{\xi_2}\xi_4 = 0.$$

Letting $\alpha = \bar{\xi_1}\xi_3$, we have $\xi_3 = \alpha\xi_1$, $\xi_4 = -\alpha\xi_2$.

Now the unitarity of the above matrix is equivalent to

$$\begin{aligned} 1 + C(\xi_1 + \xi_2) + B\eta_1 &= 0 \\ 1 + C\alpha(\xi_1 - \xi_2) - B\eta_1 &= 0 \\ B + C(\xi_1 - \xi_2)\bar{\eta_3} - \eta_1 &= 0 \\ B + C\alpha(\xi_1 + \xi_2)\bar{\eta_3} + \eta_1 &= 0. \end{aligned}$$

Thus we have

$$\xi_1 + \xi_2 = -\frac{1+\eta_1 B}{C} = -\frac{\eta_1 + B}{C}\bar{\alpha}\eta_3$$

and

$$\xi_1 - \xi_2 = \frac{-1+\eta_1 B}{C}\bar{\alpha} = \frac{\eta_1 - B}{C}\eta_3,$$

and so

$$\eta_3 = \frac{1+\eta_1 B}{\eta_1 + B}\alpha = \frac{-1+\eta_1 B}{\eta_1 - B}\bar{\alpha}.$$

This implies

$$\begin{aligned}\alpha^2 &= -\frac{(1-\eta_1 B)(\eta_1 + B)}{(1+\eta_1 B)(\eta_1 - B)} \\ &= -\frac{B^2 - 1 + B(\eta_1 - \bar{\eta_1})}{B^2 - 1 - B(\eta_1 - \bar{\eta_1})} \\ &= -\frac{2\sqrt{3} + \eta_1 - \bar{\eta_1}}{2\sqrt{3} - (\eta_1 - \bar{\eta_1})}.\end{aligned}$$

Therefore

$$\alpha = \pm i \frac{2\sqrt{3} + \eta_1 - \bar{\eta_1}}{\left|2\sqrt{3} + \eta_1 - \bar{\eta_1}\right|}.$$

□

**Corollary 3.5.** *Any biunitary connection on* 3311 *is gauge equivalent to one in the following one-parameter family of biunitary connections (using the notation from above and* $[k] = \frac{q^k - q^{-k}}{q - q^{-1}}$ *with* $[2] = \sqrt{3+\sqrt{3}}$*).*

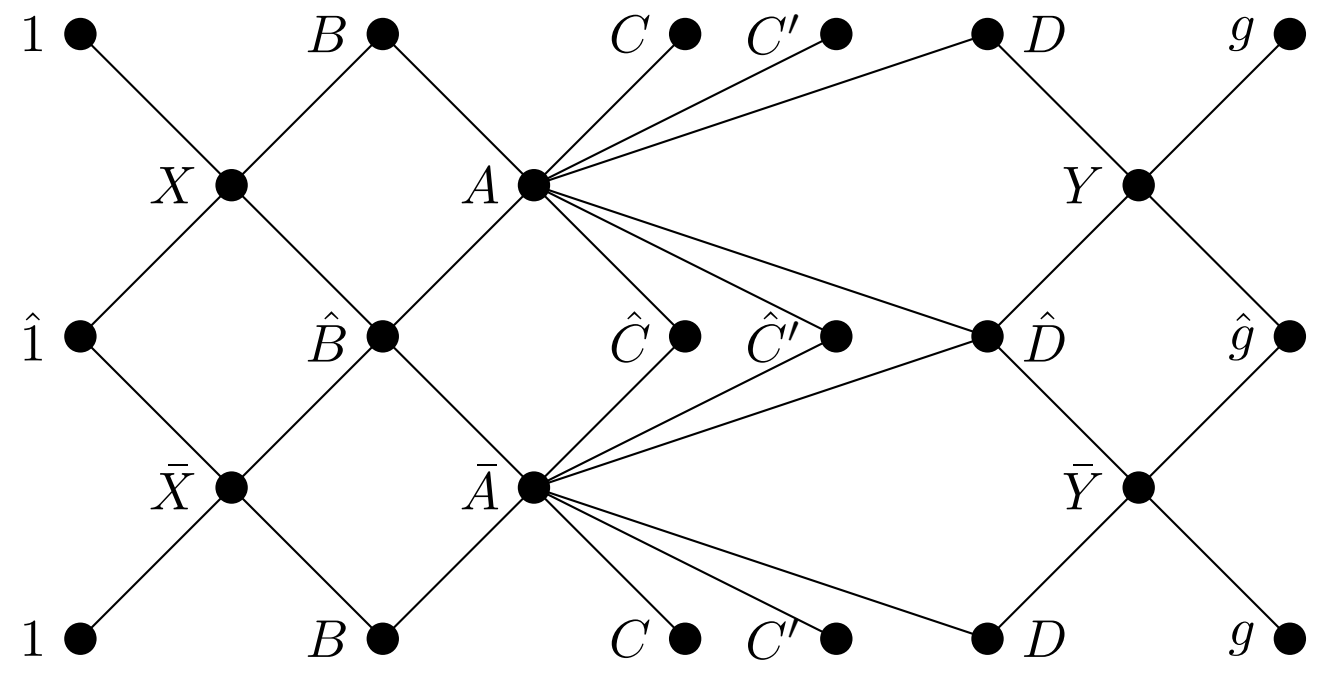

| | $1X$ | $BX$ | $BA$ | $CA$ | $C'A$ | $DA$ | $DY$ | $gY$ |
|---|---|---|---|---|---|---|---|---|
| $\hat{1}\bar{X}$ | 1 | 1 | | | | | | |
| $\hat{B}\bar{X}$ | 1 | $-\frac{1}{[3]}$ | $\frac{\sqrt{[2][4]}}{[3]}$ | | | | | |
| $\hat{B}\bar{A}$ | | $\frac{\sqrt{[2][4]}}{[3]}$ | $\frac{1}{[3]}$ | 1 | 1 | 1 | | |
| $\hat{C}\bar{A}$ | | | 1 | $\xi_1$ | $\xi_2$ | $\eta_3$ | | |
| $\hat{C}'\bar{A}$ | | | 1 | $\xi_3$ | $\xi_4$ | $-\eta_3$ | | |
| $\hat{D}\bar{A}$ | | | 1 | $\eta_1$ | $-\eta_1$ | $-\frac{1}{[3]}$ | $\frac{\sqrt{[2][4]}}{[3]}$ | |
| $\hat{D}\bar{Y}$ | | | | | | $\frac{\sqrt{[2][4]}}{[3]}$ | $\frac{1}{[3]}$ | 1 |
| $\hat{g}\bar{Y}$ | | | | | | | 1 | $-1$ |

| | $\bar{X}1$ | $\bar{X}B$ | $\bar{A}B$ | $\bar{A}C$ | $\bar{A}C'$ | $\bar{A}D$ | $\bar{Y}D$ | $\bar{Y}g$ |
|---|---|---|---|---|---|---|---|---|
| $X\hat{1}$ | $\frac{1}{[2]}$ | $\frac{\sqrt{[3]}}{[2]}$ | | | | | | |
| $X\hat{B}$ | $\frac{\sqrt{[3]}}{[2]}$ | 1 | 1 | | | | | |
| $A\hat{B}$ | | 1 | $\frac{1}{A}$ | $\frac{\sqrt{BC}}{A}$ | $\frac{\sqrt{BC}}{A}$ | $\frac{B}{A}$ | | |
| $A\hat{C}$ | | | $\frac{\sqrt{BC}}{A}$ | $\frac{C}{A}\xi_1$ | $\frac{C}{A}\xi_2$ | $\frac{\sqrt{BC}}{A}\eta_3$ | | |
| $A\hat{C}'$ | | | $\frac{\sqrt{BC}}{A}$ | $\frac{C}{A}\xi_3$ | $\frac{C}{A}\xi_4$ | $-\frac{\sqrt{BC}}{A}\eta_3$ | | |
| $A\hat{D}$ | | | $\frac{B}{A}$ | $\frac{\sqrt{BC}}{A}\eta_1$ | $-\frac{\sqrt{BC}}{A}\eta_1$ | $\frac{1}{A}$ | 1 | |
| $Y\hat{D}$ | | | | | | 1 | $\frac{1}{[2]}$ | $\frac{\sqrt{[3]}}{[2]}$ |
| $Y\hat{g}$ | | | | | | | $\frac{\sqrt{[3]}}{[2]}$ | $-\frac{1}{[2]}$ |

*Proof.* This is a straightforward application of the definition of a biunitary connection with the above lemma used to determine the 4-by-4 matrix corresponding to the quadruple point. □

In §5 below we will show that there are exactly four points on this family which give a flat connection: $\eta_1 = \pm 1$ each give the GHJ subfactor, while $\eta_1 = \pm i$ each give the dual of the GHJ subfactor. (Note that the graph automorphism negates $\eta_1$.)

## 4. Cyclotomicity

The index of a finite depth extension must be cyclotomic [dBG91, CG94, ENO05], and we next show that the index value of the $nn11$ graph (and hence, by the previous section, the index of any extension of the $(n-3)$ translate of $\mathcal{Q}^\bullet$) is never cyclotomic when $n \geq 5$. In the case of no translation we give the easy argument that the only extension of $\mathcal{Q}^\bullet$ with index $3+\sqrt{3}$ (that is, the index of 3311) is 3311 itself.

We will apply the techniques of [CMS11]. (In what follows, references to theorems and lemmas and so on are to that paper.) There are two obstacles to applying [CMS11] to $nn11$. First, we are considering a family of graphs in which we extend in two different directions at once, whereas [CMS11] considers families with a single chain of edges being extended. This obstacle is minor: the techniques of those papers would work just as well for more general families, and a graph-theoretic trick lets us ignore this subtlety entirely. The second obstacle is subtler but more significant. Although our principal graph has the same graph norm as $nn11$, there need not be any relationship between the lower eigenvalues. This leaves us in an awkward situation: Theorem 1.0.3 from [CMS11] guaranatees that the index is eventually noncyclotomic but the bounds are hopeless, while the useful bounds from Theorem 10.0.1 are irrelevant because they only guarantee that the square of some multiplicity free eigenvalue is non-cyclotomic. We will avoid this subtlety by getting lucky: the squares of the largest eigenvalue of $nn11$ turn out to be Salem numbers (see Remark 10.1.7), which means that the argument in Theorem 10.0.1 will go through.

First we explain the graph-theoretic trick. The graph $nn11$ has an involution which switches the arms of length $n$ and the arms of length 1, fixing the quadruple point. Now consider the $\mathbb{Z}/2\mathbb{Z}$ graph quotient of $nn11$,

$G_n =$ 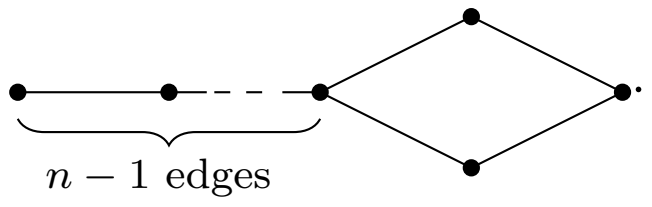

Note that $G_n$ and $nn11$ have the same graph norm (because a FP eigenvector of one gives a FP eigenvector of the other and vice-versa). Let $b_n$ be the graph norm of the $nn11$ graph (equivalently, of the $G_n$ graph).

**Theorem 4.1.** *The field $\mathbb{Q}(b_n^2)$ is not cyclotomic for any $n \geq 5$.*

*Proof.* The graphs $G_n$ are treated in [PT10], where they appear as "Vine #2", where they prove that for $n \geq 5$ the adjacency matrix of $G_n$ has a multiplicity-one eigenvalue $\lambda_n$ such that $\mathbb{Q}(\lambda_n^2)$ is non-cyclotomic. (Note that our indexing is off by 3 from the indexing there, so we need to subtract 3 from each of their bounds.) In order to prove the stronger result that $\mathbb{Q}(b_n^2)$ is non-cyclotomic, we need only observe that the characteristic polynomial of the adjacency matrix is irreducible up to cyclotomic factors. This happens because we are in what is called the "Salem case" in [PT10].

The characteristic polynomials $P_n(x)$ for the adjacency matrices of $G_n$ satisfy the recursion relation

$$P_n(x) = xP_{n-1}(x) - P_{n-2}(x)$$

with $P_3(x) = x^6 - 6x^4 + 6x^2$ and $P_4(x) = x^7 - 7x^5 + 11x^3 - 2x$. With $F_n(t) = P_n(t + t^{-1})$, we find that

$$F_n(t)\left(t - \frac{1}{t}\right) = t^n A(t) - t^{-n} A(t^{-1})$$

with $A(t) = t^4 - t^2 - 3 - 1$. This polynomial has just two real roots with magnitude greater than 1, namely $\pm\sqrt{\frac{1}{2}(1+\sqrt{17})}$, and so by Remark 10.1.7 the polynomial $P_n(x)$ is $S(x^2)$ times a product of cyclotomic polynomials where $S$ a Salem polynomial times a product of cyclotomic polynomials. Hence, if $\lambda$ is a root of $P_n$ then either $\lambda^2$ is a root of unity (and hence cyclotomic), or $\lambda^2$ is a Galois conjugate of $b_n^2$. Hence, we see that $\mathbb{Q}(b_n^2)$ must be noncyclotomic.

□

**Corollary 4.2.** *No non-trivial translate of $\mathcal{Q}^\bullet$ appears as the principal graph of a subfactor.*

*Remark.* Recall here that we always translate by an even amount, so we are only concerned here with $b_n$ for $n$ odd. As it turns out $b_4 = \sqrt{\frac{1}{2}(5+\sqrt{21})}$ is in fact cyclotomic. (Recall that the even supertransitive analogues of $\mathcal{Q}$ and $\mathcal{Q}'$ were dealt with in [MS10, Theorem 3.7] via [Jon03, Theorem 5.2.2].)

We next show that the only extension of $\mathcal{Q}^\bullet$ which can be the principal graph of a subfactor is 3311.

**Lemma 4.3.** *The only extension of $\mathcal{Q}^\bullet$ with index $3+\sqrt{3}$ is the 3311 graph.*

*Proof.* In fact, any nontrivial extension of $\mathcal{Q}^\bullet$ contains 3311, and so any extension other than 3311 itself must have index that is too large. The only other possibility is just $\mathcal{Q}^\bullet$ itself, but this does not satisfy Corollary 3.2. □

## 5. Uniqueness of the GHJ subfactor

In this section we prove the uniqueness of the 3311 subfactor up to duality. This was originally shown in unpublished work of Rehren's from 1994. For completeness we give a new proof here. Our proof combines planar algebra techniques with connections techniques and thus may be of independent interest.

The key idea of this section is that if you have a 1-parameter family of connections, then looking at rotational eigenvalues (which must be certain roots of unity) should show that only finitely many of those connections can come from a subfactor. In order to apply this technique we give a recipe which gives the rotational eigenvalues directly from the connection.

*Remark.* Uniqueness can also by proved by calculating fusion multiplicities using [KS01, Gne02]. Given that the multiplicity of $\alpha$ inside $x \otimes x$ is 0 for $R$ and 1 for $R'$, we expect that this fusion multiplicity should be expressible as a non-trivial function of $\eta_1$, and thus to give a non-trivial condition on $\eta_1$ in order for the corresponding biunitary connection to be flat. Indeed, the multiplicity of $\alpha$ inside $x \otimes x$ is 0 exactly when $\eta_1^2 = -1$, and 1 exactly when $\eta_1^2 = 1$. This gives a separate proof of the uniqueness of the 3311 subfactor up to duality. We do not include the detailed calculations here as they are somewhat lengthy.

Before moving on to the proof, we briefly discuss the structure of this subfactor from the algebraic point of view. This point of view is not necessary to understand our proof of uniqueness, but may be useful since the literature on the GHJ construction is not easily accessible to readers coming from fusion categories.

5.1. **The structure of the 3311 subfactor.** There are exactly two possible fusion rings for the even part of 3311, namely $R$ and $R'$ given in Figures 5.1 and 5.1. (Bisch [Bis94] showed that there were at most five fusion rules, but three of his rules can be eliminated because they do not satisfy $(\alpha\beta)^* = \beta^*\alpha^*$.) Each of these possible fusion rings can be realized explicitly as the even part of a subfactor as follows.

We quickly describe the GHJ construction in algebraic language. There is a restriction functor $A_{11}^{\text{even}} \to E_6^{\text{even}}$, and thus $E_6^{\text{even}}$ is a module category over $A_{11}^{\text{even}}$. Then for any simple object $\eta$ in $E_6^{\text{even}}$, we have that $\underline{\mathrm{Hom}}(\eta,\eta)$ is a $Q$-system (that is a unitary algebra object) in $A_{11}^{\text{even}}$. In particular, if we pick $\eta = 1$ we get a 3311 subfactor. (Note that although there are two different $E_6$ fusion categories, there's only one $E_6$ module category over $A_{11}$, thus doing the above construction with $\overline{E}_6$ instead gives the same subfactor.)

FIGURE 1. Fusion ring $R$

| $\otimes$ | $1$ | $B$ | $C$ | $C'$ | $D$ | $g$ |
|---|---|---|---|---|---|---|
| $1$ | $1$ | $B$ | $C$ | $C'$ | $D$ | $g$ |
| $B$ | $B$ | $1+B+C+C'+D$ | $B+C'+D$ | $B+C+D$ | $B+C+C'+D+g$ | $D$ |
| $C$ | $C$ | $B+C'+D$ | $1+C+D$ | $B+C'+g$ | $B+C+D$ | $C'$ |
| $C'$ | $C'$ | $B+C+D$ | $C'+B+g$ | $1+C+D$ | $B+C'+D$ | $C$ |
| $D$ | $D$ | $B+C+C'+D+g$ | $B+C+D$ | $B+C'+D$ | $1+B+C+C'+D$ | $B$ |
| $g$ | $g$ | $D$ | $C'$ | $C$ | $B$ | $1$ |

FIGURE 2. Fusion ring $R'$

| $\otimes$ | $1$ | $B$ | $C$ | $C'$ | $D$ | $g$ |
|---|---|---|---|---|---|---|
| $1$ | $1$ | $B$ | $C$ | $C'$ | $D$ | $g$ |
| $B$ | $B$ | $1+B+C+C'+D$ | $B+C'+D$ | $B+C+D$ | $B+C+C'+D+g$ | $D$ |
| $C$ | $C$ | $B+C'+D$ | $1+2C+g$ | $B+D$ | $B+C'+D$ | $C$ |
| $C'$ | $C'$ | $B+C+D$ | $B+D$ | $1+2C'+g$ | $B+C+D$ | $C'$ |
| $D$ | $D$ | $B+C+C'+D+g$ | $B+C'+D$ | $B+C+D$ | $1+B+C+C'+D$ | $B$ |
| $g$ | $g$ | $D$ | $C$ | $C'$ | $B$ | $1$ |

Kawahigashi [Kaw95] proved that the dual even part of the GHJ subfactor (i.e. $(A_{11}^{\text{even}})^*_{E_6^{\text{even}}}$) has fusion ring $R'$. We can also see this directly from the algebra: since the action of $A_{11}^{\text{even}}$ on $E_6^{\text{even}}$ factors through the quotient functor to $E_6^{\text{even}}$, the dual $(A_{11}^{\text{even}})^*_{E_6^{\text{even}}}$ takes a functor from $(E_6^{\text{even}})^*_{E_6^{\text{even}}} = \overline{E}_6^{\text{even}}$. Thus the dual must have the fusion rules $R'$. (Similarly, the dual takes a functor from $\overline{E}_6^{\text{even}}$. Thus the subcategories generated by $C$ and $C'$ are conjugate to each other.)

It is not difficult to show, using diagrammatic techniques, that any subfactor with fusion ring $R$ must be the example constructed above. Indeed, since $A_{11}^{\text{even}}$ comes from the quantum group SO(3), work of Tuba and Wenzl [TW05] shows that the even part of a subfactor with fusion ring $R$ must be $A_{11}^{\text{even}}$. Seeing that the odd part comes from the module category $E_6^{\text{even}}$ follows from the classification of algebra objects in $A_{11}$.

We can also give a more explicit construction of a subfactor with fusion ring $R'$. The fusion categories $E_6^{\text{even}}$ and $\overline{E}_6^{\text{even}}$ each have $\text{Vec}(\mathbb{Z}/2)$ as a fusion subcategory. In fact this fusion subcategory is central in the sense that the functor $\text{Vec}(\mathbb{Z}/2) \to E_6^{\text{even}}$ lifts to a functor $\text{Vec}(\mathbb{Z}/2) \to Z(E_6^{\text{even}})$ to the Drinfel'd center. Thus we can consider the tensor product $E_6^{\text{even}} \boxtimes_{\text{Vec}(\mathbb{Z}/2)} \overline{E}_6^{\text{even}}$. (This is the same construction that Ocneanu uses when he discusses tensoring "over the ambichirals." [Ocn01]. For more details see [Gre10, §6.1].) Note that $E_6^{\text{even}} \boxtimes \overline{E}_6^{\text{even}}$ acts on $E_6^{\text{even}}$ where the first tensor factor acts by tensoring on the left while the right tensor factor acts by tensoring on the right. Furthermore, since $\text{Vec}(\mathbb{Z}/2)$ is central (that is, comes with a fixed lift to the Drinfel'd center), this action descends to an action of $E_6^{\text{even}} \boxtimes_{\text{Vec}(\mathbb{Z}/2)} \overline{E}_6^{\text{even}}$ on $E_6^{\text{even}}$. Again, we look at $\eta$ being the trivial object in $E_6^{\text{even}}$, and the $Q$-system $\underline{\text{Hom}}(\eta, \eta)$ gives a 3311 subfactor.

However, it is not at all clear that this example is the only way to realize the fusion ring $R'$. In general the free product of two fusion categories may have many quotients which look like a tensor product on the level of objects. Furthermore, the two 3-object fusion subcategories could both be $E_6^{\text{even}}$ or both be $\overline{E}_6^{\text{even}}$. We were not able to come up with a simple algebraic argument that excludes the possibility that there is a 3311 subfactor both of whose even parts have fusion ring $R'$ and where one even part has two copies of $E_6^{\mathrm{even}}$ while the other even part has two copies of $\overline{E}_6^{\mathrm{even}}$. However, it follows from the main result of this section that no such subfactor exists.

5.2. **Rotations and connections.** To get a connection from a subfactor (or more generally, to get the $6j$-symbols in a 2-category), you first fix a collection of $3j$ intertwiners. Different choices correspond to different connections in the same gauge class. For the part of the connection attached to the initial branch point (which we call the "branch matrix"), the planar algebraic approach gives a distinguished gauge choice. On the one hand, this gauge choice is easy to recognize (it is characterized by the property that the entries in the first row and column are real numbers with specific signs). On the other hand, the diagrammatic branch matrix encodes precisely how rotation $\rho^{\frac{1}{2}}$ acts on the perpendicular complement of Temperley-Lieb.

Specifically we will show that if $U$ is the diagrammatic branch matrix, then $tr(UU^t) - 2$ is the sum of the rotational eigenvalues on the perpendicular complement of Temperley-Lieb one past the initial branch point . For 3311 these rotational eigenvalues are 1 and $-1$, so we see that $tr(UU^t) - 2 = 0$, which yields uniqueness up to duality of the 3311 subfactor.

Consider an $n$-supertransitive subfactor where the initial branch point is $k+2$-valent with no multiple edges. We will only look at the part of the connection corresponding to the initial branch point (this is a $(k+2)$-by-$(k+2)$ unitary matrix), which we will call the "branch matrix." The choice of gauge corresponds to a choice of a size 1 element in each of the 1-dimensional spaces $\mathrm{Hom}_{\pm}\big(f^{(n)} \otimes f^{(1)} \to f^{(n-1)}\big)$ and $\mathrm{Hom}_{\pm}\big(f^{(n)} \otimes f^{(1)} \to P_i^{\pm}\big)$ where the $P_i^{\pm}$ are the simple projections at depth $n+1$. Given such choices of trivalent vertices, the branch matrix $U_{\pm}$ is the change of basis matrix between the following bases of $\mathrm{Hom}_{\pm}\big(f^{(n)} \otimes f^{(1)} \to f^{(1)} \otimes f^{(n)}\big)$.

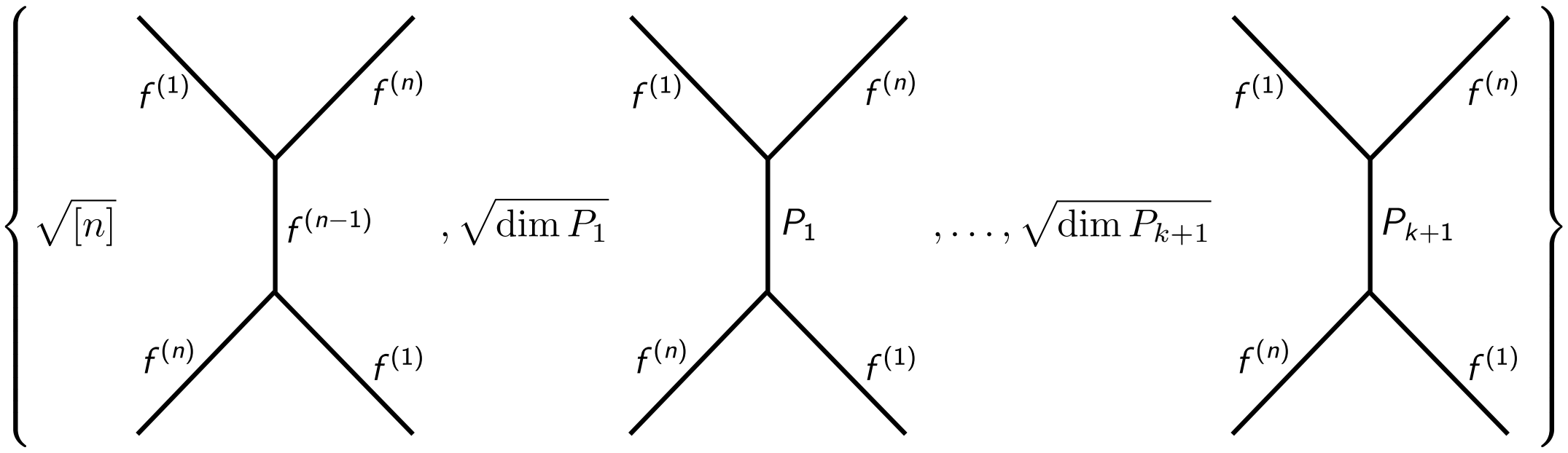

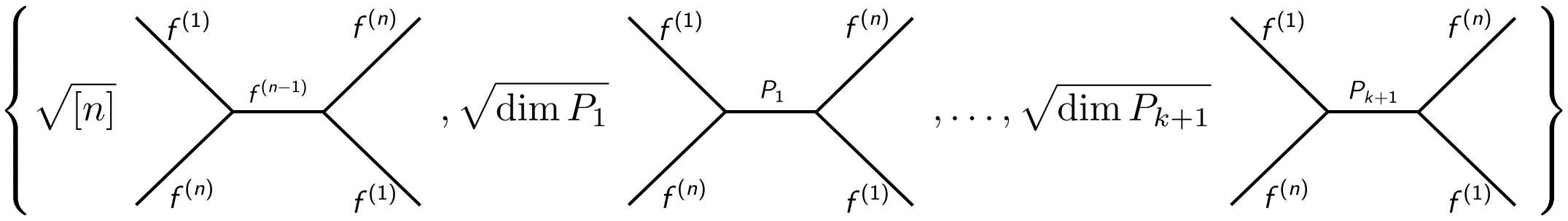

*Remark* 5.1. The multiplicative factor in the bases is there so that we can use conventions which are compatible with the usual inner product on the planar algebra. In [MPPS10] we did not have these multiplicative factors, but this was because we used a non-standard inner product.

We fix the following natural gauge choice.

**Definition 5.2.** The distinguished elements in $\mathrm{Hom}_{\pm}\big(f^{(n-1)} \otimes f^{(1)} \to f^{(n-2)}\big)$ and in $\mathrm{Hom}_{\pm}\big(f^{(n-1)} \otimes f^{(1)} \to P_i^{\pm}\big)$ are given by the following diagrams.

$\frac{1}{\sqrt{[n+1]}}$ 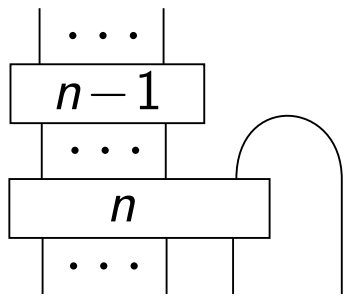 $\in \mathrm{Hom}_{\pm}\Big(f^{(n-1)} \otimes f^{(1)} \to f^{(n-2)}\Big)$

$\frac{1}{\dim P_i}$ 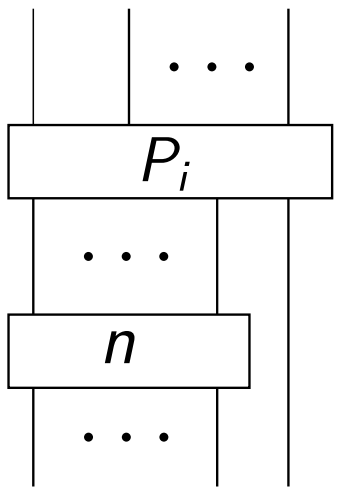 $\in \mathrm{Hom}_{\pm}\Big(f^{(n-1)} \otimes f^{(1)} \to P_i^{\pm}\Big)$

The diagrammatic branch matrix $U_{\pm}$ is the $(k+2)$-by-$(k+2)$ matrix which comes from the distinguished diagrammatic basis.

**Corollary 5.3.** *$U_+$ is the matrix for the linear transformation $\rho^{\frac{1}{2}}$ with respect to the following two bases.*

$b_1^+ =$ 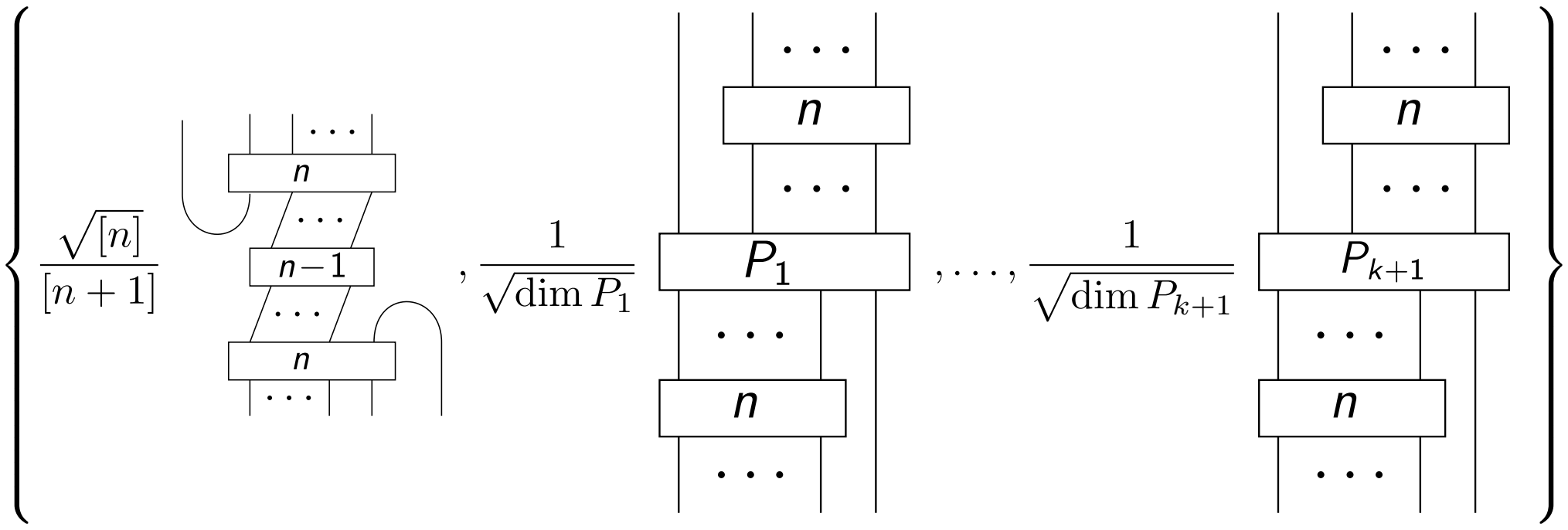

$b_2^- =$ 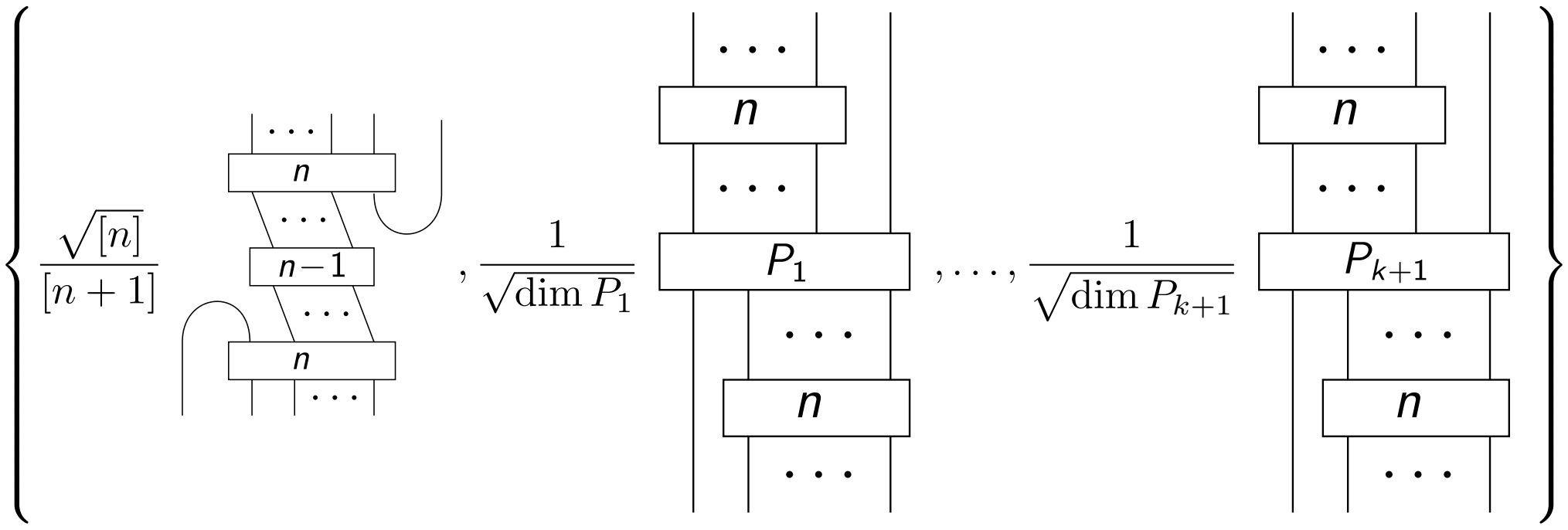

*Similarly $U_-$ is the matrix for $\rho^{\frac{1}{2}}$ going from $b_1^-$ to $b_2^+$.*

**Definition 5.4.** Let $\varphi_{\pm} : \mathbb{C}^{k+2} \to \mathcal{P}_{\pm}^{k+2}$ send the standard ordered basis of $\mathbb{C}^{k+2}$ to $b_1^{\pm}$ and $\psi_{\pm} : \mathbb{C}^{k+2} \to \mathcal{P}_{\pm}^{k+2}$ send the standard ordered basis of $\mathbb{C}^{k+2}$ to $b_2^{\pm}$. Thus the corollary above says that $U_{\pm} = \psi_{\mp}^{-1} \rho^{\frac{1}{2}} \varphi_{\pm}$.

**Lemma 5.5.** *We have the following identities in Temperley-Lieb.*

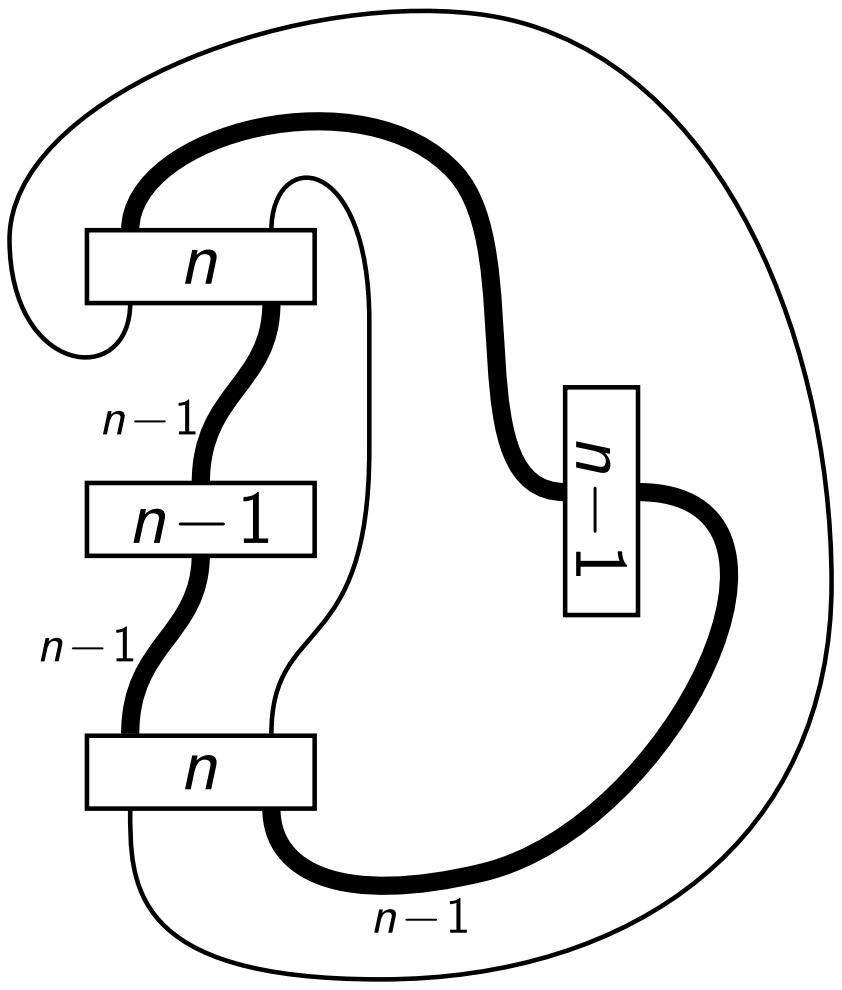

$$= (-1)^{n+1}\frac{[n+1]}{[n]}$$

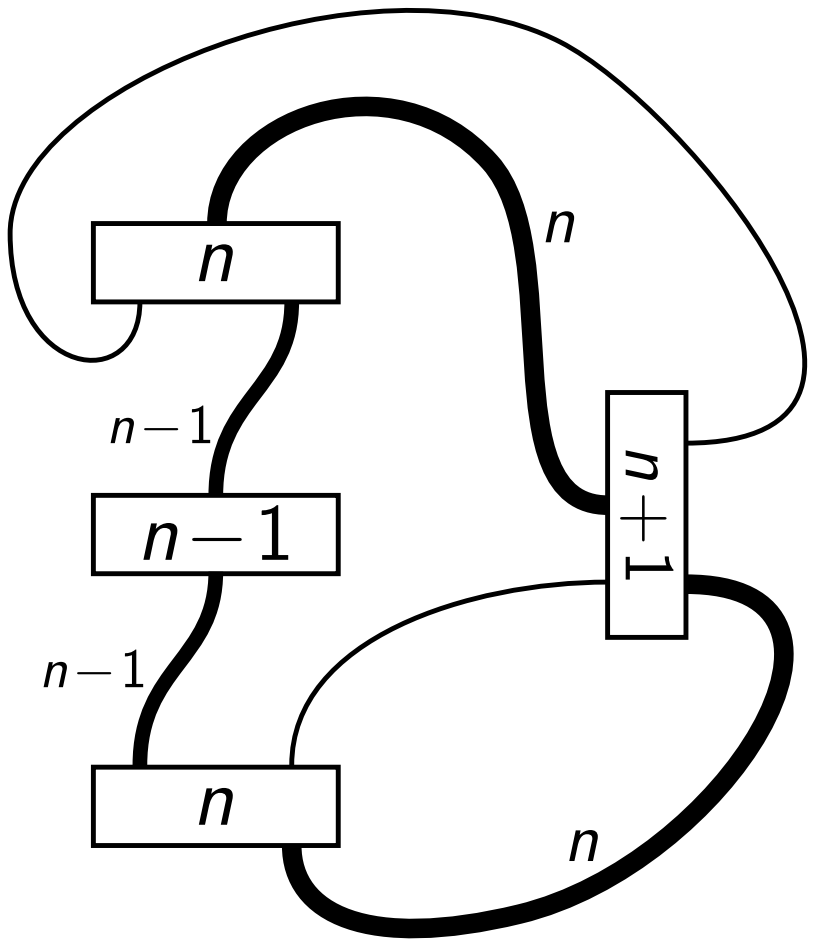

$$= [n+2]$$

*Proof.* This is a special case of the formula for tetrahedral nets in [KL94, §9.11]. (Warning: Kauffman-Lins's convention for the quantum number $[n]$ differs from the usual convention by a sign $(-1)^{n+1}$.) □

**Lemma 5.6.** *The diagrammatic branch matrix is characterized within its gauge class by the properties that the top left entry is a nonzero real number with sign $(-1)^{n+1}$, while the other entries in the first row and column are positive real numbers.*

*Proof.* First we check these properties for the diagrammatic branch matrix. The top-left entry is $\frac{\sqrt{[n]}}{[n+1]}$ times the first diagram in the previous lemma. For the other entries we need to evaluate the following diagram.

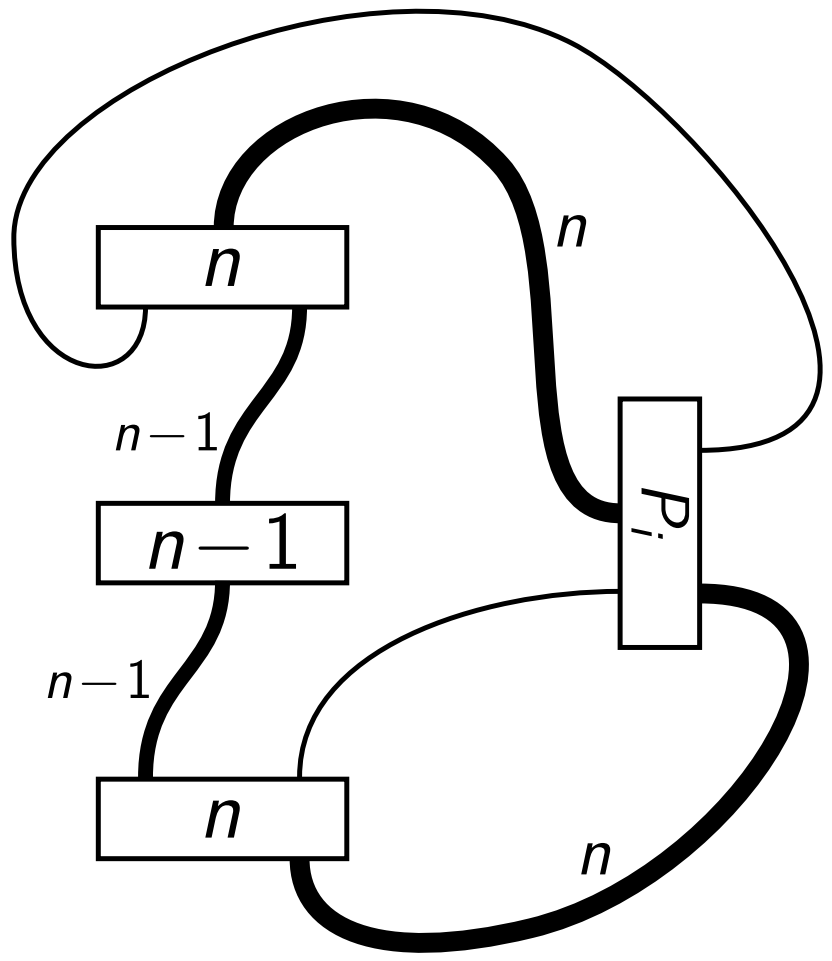

This diagram is an inner product of $P_i$ with a diagram in Temperley-Lieb. But $P_i$ is perpendicular to every Temperley-Lieb diagram except for the identity. Similarly, $f^{(n+1)}$ is perpendicular to every Temperley-Lieb diagram except for the identity. Hence,

$$\frac{\sqrt{[n]}}{[n+1]}\frac{1}{\sqrt{\dim P_i}}$$

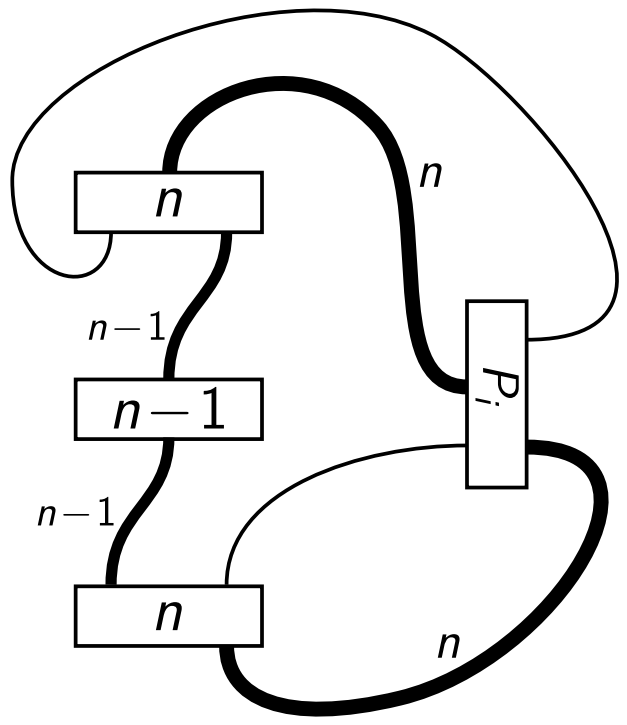

$$=$$

$$= \frac{\sqrt{[n]}}{[n+1]}\frac{1}{\sqrt{\dim P_i}}\frac{\dim P_i}{[n+2]}$$

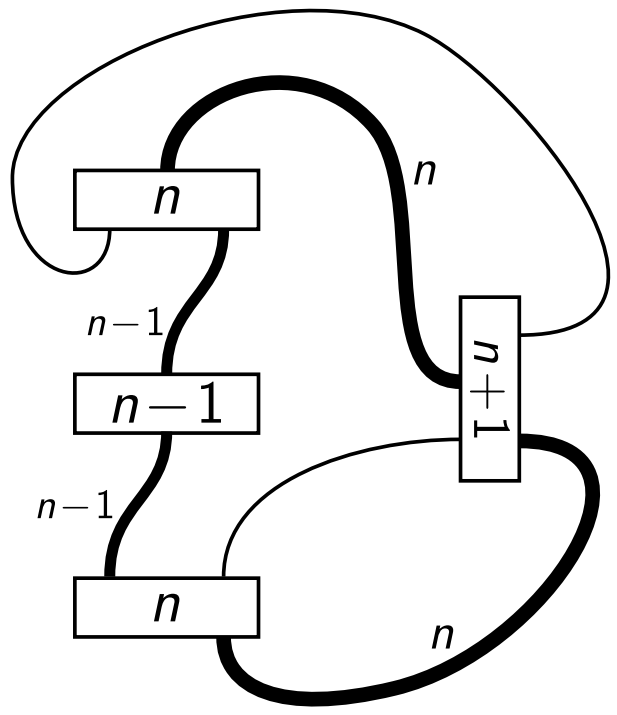

$$= \frac{\sqrt{[n]}\dim P_i}{[n+1]},$$

where the last equality follows from the value of the second diagram in the previous lemma.

Since gauge transformations rescale a row or a column by a unit complex number, exactly one branch matrix in each gauge class satisfies the conditions of the theorem. □

**Corollary 5.7.** *For each choice of $\eta_1$, the 4-by-4 matrix Lemma 3.4 is the diagrammatic branch matrix in its gauge class for a 3311 subfactor.*

**Definition 5.8.** Let $v_1^\pm = (1, 0, \ldots, 0)$ and $v_2^\pm = (0, \sqrt{\frac{[n+2]}{\dim P_1}}, \ldots, \sqrt{\frac{[n+2]}{\dim P_{k+1}}})$. Note that both $\varphi^\pm$ and $\psi^\pm$ send $v_1$ and $v_2$ to elements of Temperley-Lieb.

**Lemma 5.9.** *$U_\pm v_1^\pm = \frac{(-1)^{n+1}}{[n+1]} v_1^\mp + \frac{\sqrt{[n][n+2]}}{[n+1]} v_2^\mp$ and $U_\pm v_2^\pm = \frac{\sqrt{[n][n+2]}}{[n+1]} v_1^\mp - \frac{(-1)^{n+1}}{[n+1]} v_2^\mp$.*

*Proof.* This follows from Lemma 5.5 □

**Corollary 5.10.** *$U_\pm$ sends the perpendicular complement of $\mathbb{C}v_1^\pm \oplus \mathbb{C}v_2^\pm$ to the perpendicular complement of $\mathbb{C}v_1^\mp \oplus \mathbb{C}v_2^\mp$.*

*Proof.* This follows immediately from the unitarity of $U_\pm$. □

**Lemma 5.11.** *$\varphi^\pm$ and $\psi^\pm$ send the perpendicular complement of $\mathbb{C}v_1^\pm \oplus \mathbb{C}v_2^\pm$ to the perpendicular complement of Temperley-Lieb. Furthermore, on this subspace $\varphi^\pm = \psi^\pm$.*

*Proof.* Both $\varphi$ and $\psi$ send the perpendicular complement of $v_1$ to the span of the $P_i$. Furthermore, both $\varphi$ and $\psi$ send $v_2$ to a scalar multiple of $f^{(n+1)}$. Since the perpendicular complement of $f^{(n+1)}$ in the span of the $P_i$ is exactly the perpendicular complement of Temperley-Lieb, the first result follows. The second result follows immediately from the definition of the bases $b_1$ and $b_2$. □

**Lemma 5.12.** *$U_- = U_+^t$.*

*Proof.* This follows from the normalization condition on a biunitary connection. □

Henceforth we will let $U = U_+$. We now recall Theorem 1.7 from the introduction, and provide its proof.

**Theorem 1.7.** *Suppose $U$ is the diagrammatic branch matrix for an $n$-supertransitive principal graph pair $\Gamma$, with the initial branch point is $k+2$-valent, with no multiple edges, and a flat biunitary connection.*

*The eigenvalues for $U^tU$ with multiplicity are: two 1's together with the eigenvalues for $\rho$ acting on the perpendicular complement of Temperley-Lieb in the $n$-box space. In particular, $tr(U^tU) - 2$ is the sum of the eigenvalues of $\rho$.*

*Proof.* An explicit calculation using Lemmas 5.9 and 5.12, shows that $v_1$ and $v_2$ are eigenvalue 1 eigenvectors for $U^tU$.

On the perpendicular complement of Temperley-Lieb we have:

$$\rho = \rho^{\frac12}\rho^{\frac12} = \psi_+ U^t \varphi_-^{-1} \psi_- U \varphi_+^{-1} = \varphi_+ U^t \psi_-^{-1} \psi_- U \varphi_+^{-1} = \varphi_+ U^t U \varphi_+^{-1}$$

Here the third equality follows from Lemma 5.11 together with the fact that $\rho^{\frac12}$ sends $\mathcal{TL}_+^\perp$ to $\mathcal{TL}_-^\perp$.

Thus, on the perpendicular complement of $v_1$ and $v_2$, we have that $U$ has the same eigenvalues as $\rho$ on the perpendicular complement of Temperley-Lieb. □

*Remark.* Applying the above argument to the dual subfactor we also see that $UU^t$ has the same properties.

*Example.* Using the notation and results of [AH99], the diagrammatic branch matrix for the dual Haagerup subfactor is :

$$U = \begin{pmatrix} \frac{1}{\lambda(\lambda^2-2)} & \frac{1}{\sqrt{3}} & \frac{\sqrt{\lambda^2-1}}{\lambda^2-2} \\ \frac{\lambda^2-1}{\lambda(\lambda^2-2)} & \frac{1}{\sqrt{3}}\tau\bar{\rho} & \frac{1}{\lambda^2-2}\bar{\rho} \\ \frac{\lambda^2-1}{\lambda(\lambda^2-2)} & \frac{1}{\sqrt{3}}\tau\rho & \frac{1}{\lambda^2-2}\rho \end{pmatrix}$$

$$= \begin{pmatrix} \frac{\sqrt{4-\sqrt{13}}}{3} & \frac{1}{\sqrt{3}} & \frac{\sqrt{2+\sqrt{13}}}{3} \\ \frac13\sqrt{\frac{5+\sqrt{13}}{2}} & -\frac12\sqrt{\frac13\left(7-3\sqrt{13}-i\sqrt{6\left(-25+7\sqrt{13}\right)}\right)} & -\sqrt{-\frac{31}{36}+\frac{7\sqrt{13}}{36}+\frac12 i\sqrt{\frac16\left(-25+7\sqrt{13}\right)}} \\ \frac13\sqrt{\frac{5+\sqrt{13}}{2}} & -\frac12\sqrt{\frac13\left(7-3\sqrt{13}+i\sqrt{6\left(-25+7\sqrt{13}\right)}\right)} & -\sqrt{-\frac{31}{36}+\frac{7\sqrt{13}}{36}-\frac12 i\sqrt{\frac16\left(-25+7\sqrt{13}\right)}} \end{pmatrix}$$

A direct computation shows that

$$U^tU = \begin{pmatrix} 1 & 0 & 0 \\ 0 & \frac{3-\sqrt{13}}{2} & \sqrt{\frac{-9+3\sqrt{13}}{2}} \\ 0 & \sqrt{\frac{-9+3\sqrt{13}}{2}} & \frac{-3+\sqrt{13}}{2} \end{pmatrix}$$

and

$$UU^t = \begin{pmatrix} 1 & 0 & 0 \\ 0 & 0 & 1 \\ 0 & 1 & 0 \end{pmatrix}.$$

An eigenbasis of $U^tU$ consists of the two Jones-Wenzl projections $(\sqrt{\frac{3+\sqrt{13}}{2}}, 0, 0)$ and $(0, \sqrt{\frac{1+\sqrt{13}}{2}}, \sqrt{\frac{5+\sqrt{13}}{2}})$ with eigenvalue 1, and the rotational eigenvector $(0, \sqrt{\frac{5+\sqrt{13}}{2}}, -\sqrt{\frac{1+\sqrt{13}}{2}})$ with eigenvalue $-1$. Recall that $(0, \sqrt{\frac{5+\sqrt{13}}{2}}, -\sqrt{\frac{1+\sqrt{13}}{2}})$ corresponds to

$$\sqrt{\frac{5+\sqrt{13}}{2}}\frac{Q'}{\sqrt{\dim Q'}} - \sqrt{\frac{1+\sqrt{13}}{2}}\frac{P'}{\sqrt{\dim P'}} = \sqrt{r}Q' - \frac{1}{\sqrt{r}}P',$$

where $r = \frac{\dim P'}{\dim Q'}$. This is proportional to the formula for the rotational eigenvector in [Jon03] where it was also shown that the rotational eigenvalue is $-1$.

Similarly, an eigenbasis of $UU^t$ consists of the two Jones-Wenzl projections $(\sqrt{\frac{3+\sqrt{13}}{2}}, 0, 0)$ and $(0, \sqrt{\frac{3+\sqrt{13}}{2}}, \sqrt{\frac{3+\sqrt{13}}{2}})$ with eigenvalue 1, and the rotational eigenvector $(0, \sqrt{\frac{3+\sqrt{13}}{2}}, -\sqrt{\frac{3+\sqrt{13}}{2}})$ with eigenvalue $-1$. In this case, the rotational eigenvector is the element $P - Q$.

5.3. **Uniqueness.**

**Lemma 5.13.** *For any 3311 subfactor planar algebra, the eigenvalues of $\rho$ acting on the perpendicular complement of Temperley-Lieb in the $n$-box space are 1 and $-1$.*

*Proof.* The eigenvalues are 4th roots of unity. By [Jon03, Theorem 5.2.3.], the eigenvalues are distinct. Since $C$ and $C'$ are self-dual (either by Theorem 1.3 or by looking at the fusion rings $R$ and $R'$), neither eigenvalue can be $\pm i$. Hence they must be 1 and $-1$. □

**Corollary 5.14.** *$tr(UU^t) = 1 + 1 + 1 - 1 = 2$.*

**Lemma 5.15.** *With $U$ the four-by-four unitary matrix from Lemma 3.4 depending on the modulus 1 parameter $\eta_1$ and the choice of sign for $\alpha$, we have that $tr(UU^t) = 2$ if and only if $\eta_1^4 = 1$.*

*Proof.* First, with $U$ the four-by-four matrix from Lemma 3.1 we have

$$\mathrm{tr}(UU^t) - 2 = \frac{1}{A^2}\left(x^2 + y^2 + 2B^2 + 4BC - 2A^2 + BC(\eta_1^2 + \eta_2^2 + \eta_3^2 + \eta_4^2) + C^2(\xi_1^2 + \xi_2^2 + \xi_3^2 + \xi_4^2)\right).$$

For 3311 we have $A = (1+\sqrt{3})\sqrt{3+\sqrt{3}}$, $B = 2+\sqrt{3}$, $C = 1+\sqrt{3}$ and $y = 1$, $x = -1$, and this simplifies to

$$\mathrm{tr}(UU^t) - 2 = \frac{1}{3+\sqrt{3}}\left(\frac{1+\sqrt{3}}{2}(\eta_1^2 + \eta_2^2 + \eta_3^2 + \eta_4^2) + (\xi_1^2 + \xi_2^2 + \xi_3^2 + \xi_4^2)\right).$$

Substituting in the formulas for $\eta_2$, $\eta_3$, $\eta_4$ and the $\xi_i$ from Lemma 3.4, and performing some lengthy algebraic simplifications, we obtain

$$\mathrm{tr}(UU^t) - 2 = \frac{1+\sqrt{3}}{2}\frac{\eta_1^4 - 1}{\eta_1^2 - (7+4\sqrt{3})}.$$

The denominator here is bounded away from zero, so we conclude that $\mathrm{tr}(U^tU) = 2$ exactly when $\eta_1^4 = 1$. □

**Corollary 5.16.** *If the connection from Corollary 3.5 is flat then* $\eta_1^4 = 1$.

*Proof of Theorem 1.4.* For each choice of $\eta_1$ there are exactly two connections (depending on the choice of sign for $\alpha$). However, the graph automorphism of 3311 acting on the dual graph interchanges the middle two columns of the branch matrix, this says that the choice of sign does not change the subfactor. Furthermore, the graph automorphism of 3311 acting on the principal graph negates $\eta_1$, so the subfactor only depends on the sign of $\eta_1$. Thus there are at most two subfactors with principal graph 3311: the one whose connection has $\eta_1 = \pm 1$ and the one whose connection has $\eta_1 = \pm i$.

It is easy to see that if $\eta_1 = \pm 1$, then $\eta_3 = \pm i$ and vice-versa. Thus these two subfactors (if they exist) are dual to each other. Thus up to duals there is at most one 3311 subfactor and we know that this subfactor is not self-dual. This proves uniqueness. □

*Remark.* It is not hard to see which connection corresponds to which fusion category. The even part corresponding to $\eta_1 = \pm 1$ is isomorphic to its complex conjugate under an isomorphism which fixes the two middle rows, while the even part corresponding to $\eta_3 = \pm 1$ is isomorphic to its complex conjugate under an isomorphism which interchanges the two middle rows. Thus the former case corresponds to $A_{11}^{\text{even}}$ while the latter case corresponds to $E_6^{\text{even}} \boxtimes_{\text{Vec}(\mathbb{Z}/2)} \overline{E}_6^{\text{even}}$.

This paper is available online at the arXiv, and at http://tqft.net/index5-part3.

*E-mail address*: izumi@math.kyoto-u.ac.jp

*E-mail address*: vfr@math.berkeley.edu

*E-mail address*: scott@tqft.net

*E-mail address*: nsnyder@math.columbia.edu